\documentclass[11pt]{article}
\usepackage{amsfonts}
\usepackage{amsmath}
\usepackage{amsxtra}
\usepackage{amsthm}
\usepackage{euscript}
\usepackage{amssymb}
\usepackage[T2A]{fontenc}

\newtheorem{Theorem}{Theorem}

\theoremstyle{definition}

\theoremstyle{remark}

\newcommand{\emk}{\mathop{\rm cap}\nolimits}

\newcommand{\diam}{\mathop{\rm diam}\nolimits}
\newcommand{\dist}{\mathop{\rm dist}\nolimits}

\title{
\hbox to \textwidth{\hfill}
Behavior of bounded solutions of quasilinear elliptic equations on Riemannian manifolds}
\author{A.B.Ivanov}
\date{}

\begin{document}

\large

\maketitle
				Let $M$ be a connected complete noncompact Rimannian manifold of dimension $n$, $K$ is a compact subset of $M$. We study bounded in $M\setminus K$ solutions of the equation

				\begin{equation}
				        \label{MainEquation}
									\Delta_{p} u = c(x){\left| u \right|}^{p-2} u,
				\end{equation}

				where $c(x)\in L_{\infty, loc}\left(M\right)$ is a nonnegative function, $p>1$ is a real number, and operator $\Delta_{p}$ is defined by
				
        $$
        	\Delta_{p} u =  {g}^{-\frac{1}{2}} \nabla_{i}\left( g^{ij} \left| \nabla u \right|^{p-2} \nabla_j u \right)
        $$
        where $g^{ij}$ are contravariant components of the metric tensor, $g={\rm det}\left\|g^{ij}\right\|$, and
        $$
        	\left| \nabla u \right| = \left(  g^{ij} \nabla_{i}u \nabla_{j}u\right)^{\frac{1}{2}}.
        $$
				
				Similar problems were studied in \cite{KL}--\cite{GR}.
				
				A function $u \in {W_{p, loc}^1}({\Omega}) \cap {L_{\infty}({\Omega})}$ is said to be a solution of $(\ref{MainEquation})$ in $\Omega$ if
        
				$$
					- \int\limits_{{\Omega}}{ 
																										{g^{ij} \left| \nabla u \right|}^{p-2}
																										\nabla_{i}u
																										\nabla_{j}\varphi
																										} \, dx = 
					\int\limits_{{\Omega}}{ 
																										c(x) {\left| u \right|}^{p-2} u \varphi
																										} \, dx
				$$
				for every $\varphi \in C_{0}^{\infty} ({\Omega})$.

				We say a function $u \in {W_{p, loc}^1}(\Omega) \cap {L_{\infty}(\Omega)}$ to be \mbox{$p$ -harmonic} in $\Omega$, if 
				$$
						\Delta_{p} u = 0
				$$
				in the above sense.

				Riemannian manifold $M$ is said to be an $H_p$-manifold if, for all non-empty open sets \mbox{$\omega\subset\Omega\subset M$} such that
				$$
					{\gamma}_1\le\frac{\diam(\omega)}{\diam(\Omega)}\le{\gamma}_2
				$$ 
				and
				$$
					{\gamma}_1\le\frac{\dist(\omega, \partial\Omega)}{\diam(\omega)}\le{\gamma}_2,
				$$
				where $\diam(\omega)$, $\diam(\Omega)$ are the diameters of the sets $\omega$ è $\Omega$ correspondingly, $\dist(\omega, \partial\Omega)$ is the distance between 		$\omega$ and $\partial\Omega$,  ${\gamma}_1$, ${\gamma}_2 > 0$ are real numbers,
 				and for every nonnegative function $u$ that is p-harmonic on $\Omega$, the following inequality holds:
 				$$\sup_{\omega} u \le h \inf_{\omega} u,$$
 				where the constant $h > 0$ depends only on $\gamma_1$ and $\gamma_2$.

				Let $A$ be a compact, and $B$ an open subsets of $M$, $A \subset B$. The \mbox{$p$ - capacity} of $A$ with respect to $B$ is defined by
				$$
					\emk(A,B)=\inf_{\varphi}\int\limits_{B} {\left| \nabla\varphi \right|}^{p} \, dx,
				$$
				where the infinum is taken over all functions $\varphi \in C_{0}^{\infty} (B)$ identically equal to 1 in a neighborhood of $A$.

				Fix an arbitrary point $O\in M$.

        We denote by $B_k = \{ x : \dist(O, x) < 2^k \}$ the geodesic ball of radius $2^k$ with center at $O$. For every set $\Omega$ denote 
					$Q_k=\overline{B_{k}}\setminus\Omega$,
					$D_k={B_{k+2}\setminus \overline{B_{k+1}}}$ and $d_k={\overline{B_{k+3}}\setminus B_{k}}$.

				An open set $\Omega$ is said to be \mbox{$p$ - massive} if 
				$$
					\sum_{k=0}^{\infty} {
												{\left(\frac{\emk\left(Q_{2k}, B_{2k+2}\right)}{\emk\left(\overline{B_{2k+1}}, B_{2k+2}\right)}\right)}^{\frac{1}{p-1}}
															} < \infty.
				$$

				The main result of this paper is the following theorem.
 
				\begin{Theorem} \label{T1}
				Let $M$ be a {\mbox{$H_p$ -manifold}} and $c(x)$ satisfy
				$$
								\sum_{k=0}^{\infty} {
									\left(
												\int\limits_{d_k} {\frac{c(x)}{\emk\left(d_k, D_k\right)}}  \, dx\right)^{\frac{1}{p-1}}											
								}= \infty
				$$ 
				for any \mbox{$p$ - massive} set $\Omega$.
				Then any bounded solution of $(\ref{MainEquation})$ in $M\setminus K$ tends to zero when $\dist(O,x)\to\infty$
				\end{Theorem}

\end{document}